\documentclass[11pt]{article}
\usepackage{epsfig,amssymb}

\begin{document}

\title{
Projective invariants of linear 3-webs and Gronwall's Conjecture}

\author{{\Large Sergey I. Agafonov}\\
Department of Mathematics,\\
S\~ao Paulo State University-UNESP,\\ S\~ao Jos\'e do Rio Preto, Brazil\\
e-mail: {\tt sergey.agafonov@gmail.com} }
\date{}
\maketitle
\unitlength=1mm

\newtheorem{theorem}{Theorem}
\newtheorem{proposition}{Proposition}
\newtheorem{lemma}{Lemma}
\newtheorem{corollary}{Corollary}
\newtheorem{definition}{Definition}
\newtheorem{example}{Example}

\pagestyle{plain}

\begin{abstract}
\noindent We present a projectively invariant description  of planar linear 3-webs. For a non-hexagonal 3-web, we introduce family of projective torsion-free Cartan connections, the web leaves being geodesics for each member of the family, and give a web linearization criterion. Finally, we propose an algorithm for resolving the  Gronwall conjecture and illustrate this approach by proving the conjecture for 3-webs whose 2 foliations are  2 pencils of lines.\\

\noindent {\bf Key words:} linear 3-web, Gronwall conjecture, Cartan connection. \\
\\
{\bf AMS Subject classification:} 53A60
\end{abstract}

\section{Introduction}
This paper is devoted to one of the oldest  topics of the web theory, namely, to the problem of web linearization. 

A planar 3-web $\mathcal{W}_3$  is a superposition of three foliations in the plane.
If the leaves of all three foliations are rectilinear then the web is called linear. 
A linearization of a planer 3-web $\mathcal{W}_3$ is a local diffeomorphism mapping $\mathcal{W}_3$ to some linear web $\mathcal{L}_3$. A web is called
hexagonal (or flat) if it admits a linearization sending the leaves  
 of each foliation to parallel lines. Since projective transformations  map a linear 3-web to a linear 3-web, linearization, if there is any, is not unique. Two linearizations $\varphi,\psi$ of a 3-web are projectively equivalent if there exists a projective transformation  $G\in PGL(3)$ such that $\psi=G\circ \varphi$. In what follows, a planar 3-web will be called {\it polymorphic} if it admits at least two projectively non-equivalent linearizations. 

Each foliation of a linear 3-web determines a curve arc in the dual plane.
Graf and Sauer (see \cite{GSg}) gave the following complete (and very elegant!) description of linear hexagonal 3-webs: a linear web $\mathcal{L}_3$ is hexagonal if and only if the three arcs, corresponding to the foliations,  belong to some (possibly singular) cubic. 

Thus, the lines of a hexagonal linear 3-web envelop an algebraic curve of 3d class.  
Such curves, being dual to cubics, have projective moduli. On the other hand, any two hexagonal 3-webs are locally diffeomorphic.

Gronwall conjectured in 1912 (see \cite{Gn}) that any polymorphic 3-web is hexagonal. Or, in its original form, the conjecture claims that for a non-hexagonal planar 3-web there exists at most one projective  class of linearizations.  (It is noteworthy that in this very paper Gronwall promised to prove his claim in a subsequent paper.)

Since it is believed that the conjecture is true, most of the efforts were concentrated on proving it or on finding an upper bound for the number of projectively non-equivalent linearizations 
  (see \cite{BB,Be,Bk,Bg,Bp,GLb,GMS,Vn}). Bol  \cite{Bk} found the first estimate of 17 for projective
 linearization classes,  Bor\r{u}vka \cite{Bp} lowered it to 16, the short note \cite{Vn} of
 Vaona  presents a sketch of proof that the bound is 11. 
 
 Remarkably, G.Bol  gave the following motivation for publishing his paper \cite{Bg}:
{\it "In dieser Note m\"ochte ich das angedeutete Rechenverfahren sowie die
Beispiele bekanntgeben, vor allem in der Hoffnung, dass es einem besseren
Rechner gelingen k\"onnte, hieran anschliessend die vermutete Eindeutigkeit
durch ein Gegenbeispiel zu widerlegen.\footnote{In this note, I will present the computational approach mentioned above as well as a few examples in hope that, with their help, a better calculator would manage to refute the conjectured uniqueness by a counterexample.}}"  
 
Most of the bounds were obtained as a by-product of searching for a linearizability criterion. 
Gronwall himself published the first such criterion in \cite{Gn}. The existence of linearization was reduced to the existence of solution to some weakly overdetermined non-linear system of partial differential equations (PDEs), the solution being a complete projective invariant of a searched-for linearization.   
Since each author has its own taste in choosing this invariant, the criterion came up in many guises by various authors, but, in its essence, it remains the same: PDEs for the multi-dimensional Schwarzian derivative (see \cite{Sp} for the definition of multi-dimensional Schwarzian). The main difficulty in applying these criteria lies in the nature of the obtained PDE system: being weakly overdetermined, it needs several prolongations; being non-linear, it  leads very quickly to huge polynomial compatibility conditions, intractable  even by modern computer algebra software. 

Algorithmic approach to the linearization problem was suggested in \cite{GMS} and \cite{GLb}: the idea was to go through the compatibility analysis and to find a linearizabilty criterion in terms of differential invariants of the web, i.e. to exclude Schwarzian derivative components from the PDEs of the criterion. The output do not seem very satisfactory since the final formulas are immensely involved, moreover, the results of the two mentioned papers do not match: the authors of \cite{GLb} claim that the main example of \cite{GMS} is not linearizable, whereas this example possesses an explicit linearization in elementary functions!\footnote{I thank J.P.Dufour for communicating this explicit linearization.} Anyway, both papers agree that  the bound for projective classes of linearizations is 15.   

Gronwall's conjecture was proven for some restrictions  on the map and/or on the web.

Bol \cite{Bg} showed that:\\
1) a local diffeomorphism, mapping a pencil of lines into a pencil of lines and preserving linearity of some two other foliations, is projective,\\
2)  a local diffeomorphism, mapping a linear 3-web, whose two families of  lines are tangents of some conic, to some 3-web of the same type, is projective,\\
3) a local diffeomorphism, mapping a linear 3-web, whose two families of  lines are tangents of some conic and the 3d family is arbitrary, to some linear 3-web so that the image of the 3d family is a pencil of lines, is projective.

Wang \cite{Wg} demonstrated that a polymorphic 3-web is hexagonal,
provided that its Blaschke curvature vanishes to order three at some point.

In a short note \cite{Su}, Smirnov rediscovered the Bol result 1) mentioned above, and suggested 
a line of attack on the general case, claiming that it is reducible to a web with a pencil of lines.\footnote{The reviewers of both Mathematical Reviews and Zentralblatt MATH erroneously accepted this very non-precise program as a proof of the conjecture.} 
Later Smirnov published a more detailed paper \cite{Sa} proving the Bol result 1), but never returned to his project outlined in \cite{Su}.

Finally, Gronwall's conjecture turned out to be true  for webs admitting an infinitesimal symmetry (see \cite{Agr}).

We start this paper by constructing projective differential invariants for a linear planar 3-web. The approach is classical: we adjust a projective frame to the web and recover a complete invariant as the Darboux derivative. The Darboux derivative  satisfies $SL(3)$-structure equations, thus giving differential equations for the invariants.

Then, following the tradition, we present a version of the linearizability criterion. Using the form of the obtained Darboux derivative as a model, we introduce a family of projective torsion-free Cartan connections, parametrized by one non-vanishing function, the web leaves being geodesics for each member of the family. Then the web is linearizable if and only if there is a choice of the functional parameter that kills the curvature of the connection.  The zero curvature condition impose two PDEs on the parameter. It is known that the compatibility analysis of these PDEs leads to polynomial constraints. Therefore any non-hexagonal planar 3-webs carries a finite number of "natural" projective connections, which are candidates for verifying the linearization criterion. 

In the 9-dimensional space of projective differential invariants of second order,  the invariants of a particular 3-web parametrize some 2-dimensional surface, which we call {\it the signature set} of the web. Two linear 3-webs are projectively equivalent if and only if they have the same signature set. The signature set can degenerate to a point or to a curve. The former degeneration corresponds to  hexagonal linear 3-webs with 2-dimensional projective symmetry, the latter to linear 3-webs with 1-dimensional projective symmetry.

Further we derive differential equations for maps  preserving linearity of a given linear 3-web, and analyze their compatibility conditions. The novelty of our approach is in considering the obtained system as defining a polymorphic web: we do not try to exclude  projective invariants  of the map (i.e. its Schwarzian derivative in a new guise) in the first place. This allows us to go a bit further through the compatibility analysis  and to estimate, for example, the dimension of moduli space for polymorphic 3-webs. Moreover, the analysis shows that polymorphic webs are necessarily analytic.

Finally, we propose an algorithm that, being implemented on a sufficiently powerful computer, will resolve the Gronwall conjecture. The algorithm is based on two facts:
 
 1) the conjecture is true for webs with infinitesimal symmetry \cite{Agr}, 

 2) a web with one-dimensional signature set admits an infinitesimal symmetry.\\
  We illustrate this approach applying it to 3-webs whose two foliations are pencils of lines; the conjecture turns out to be true for such webs.  

As for the general case, it seems to be out of reach for computer facilities available. Noteworthily, pioneers in the field were quite right in their estimates of the computational difficulty of the problem, namely, G.Bol and W.Blaschke left the following comment  in \cite{BB}: {\it "Dieses Eindeutigkeitsproblem ... wird sich wohl kaum l\"osen lassen, solange nicht in den Vereinigten Staaten die entsprechende Rechenmaschine konstruirt worden ist.\footnote{This uniquiness problem will hardly be resolved until an appropriate computer is constructed in the United States.}"}

The interest towards the Gronwall conjecture is explained by its importance for the projective geometry foundations: namely, if the conjecture is true then the topological structure of a linearizable planar non-hexagonal 3-web determines a unique projective structure, at least locally.               
  Note that this is true for $4$-webs. Any $4$-web supplies the underlying manifold with  a unique
projective connection such that the web leaves are geodesic.  Linearizability of the 4-web amounts to the flatness of the corresponding connection
(see \cite{Lg,BB}) for the classical and \cite{Hl,AGL,Pl} for a modern treatment).

All the objects treated in this paper are local and smooth, the results are true in both real and complex settings. 

\section{Construction of differential invariants via Darboux derivative}

Let a planar linear 3-web $\mathcal{L}_3$ be defined on some open set $V$ of the projective plane $\mathbb P^2$. Each foliation $\mathcal{F}_i$ of the web, being a family of straight lines, envelops some focal curve $C_i$, which may degenerate to a point. Each point $p \in V$ belongs to some line $l_i(p)\in \mathcal{F}_i$ of the $i$-th family, the line touching the corresponding focal curve $C_i$ at some well-defined point $\eta_i(p)\in C_i$. Observe that for  webs, defined only locally, the focal curves need not belong to $V$.
\begin{lemma}\label{aligned}
If the 3 points $\eta_1(p),\eta_2(p),\eta_3(p)$ are collinear for any point $p\in V$ then all focal curves $C_i$ degenerate to points.
\end{lemma}
 {\it Proof:} Let us choose an affine chart containing the points $\eta_1(p),\eta_2(p),\eta_3(p)$ and denote by $P,Q$ and $R$ the inclinations of the 3 web lines, meeting at the point with the affine coordinates $(x,y)$. Then each of the direction fields $P(x,y),Q(x,y),R(x,y)$ satisfies the Euler equation:
 \begin{equation}\label{Euler}
P_x+PP_y=0,\ \ \
Q_x+QQ_y=0,\ \ \
R_x+RR_y=0.
\end{equation}
In the chosen affine coordinates $x,y$, one finds
$$
\begin{array}{c}
\eta_1(x,y)=(x-\frac{1}{P_y},y-\frac{P}{P_y}),\\
\eta_2(x,y)=(x-\frac{1}{Q_y},y-\frac{Q}{Q_y}),\\
\eta_3(x,y)=(x-\frac{1}{R_y},y-\frac{R}{R_y}).\\
\end{array}
$$
Collinearity of these points amounts to
\begin{equation}\label{constraintLine}
P_y(Q-R)+Q_y(R-P)+R_y(P-Q)=0.
\end{equation}
Computing the compatibility conditions of this differential constraint with system (\ref{Euler}) one gets $P_{yy}=Q_{yy}=R_{yy}=0$. Now from  (\ref{Euler}) we obtain  $P(x,y)=\frac{y-y_1}{x-x_1}$, $Q(x,y)=\frac{y-y_2}{x-x_2}$, $R(x,y)=\frac{y-y_3}{x-x_3}$. Therefore $\eta_i(x,y)=(x_i,y_i)$ and
all $\eta_i(p)$ are stable.
The details of  computing the compatibility conditions are presented in the Appendix.
\hfill $\Box$\\

\noindent{\bf Remark.} Actually, above we have proved also that the web is formed by 3 pencils of straight lines if and only if $P_{yy}\equiv Q_{yy}\equiv R_{yy}\equiv 0$. If the pencil centers are collinear then the web is called {\it regular}.\\

\noindent Choose some $\zeta(p), \ \xi_i(p) \in \mathbb K^3$, where $\mathbb  K=\mathbb C$ or $\mathbb K=\mathbb R$, to represent the points $p$ and $\eta_i(p)\in C_i$ and denote $F(p):=(\xi_1(p),\xi_2(p),\xi_3(p))$ the matrix composed of vector-columns $\xi_i(p)$. As follows from Lemma \ref{aligned}, for non-regular 3-web,  the vectors $\xi_1(p),\xi_2(p),\xi_3(p)$ form a basis in $\mathbb K^3$ and one can always normalize them so that
\begin{equation}\label{norm}
 \zeta(p) \wedge[\xi_1(p)+\xi_2(p)+\xi_3(p)]=0, \ \ \ \det F(p)\equiv 1.
\end{equation}
Explicitly, one computes  
$$
F(p)=\frac{1}{\sqrt[3]{\mu}}\left(
  \begin{array}{ccc}
    z_1 \left(x-\frac{1}{P_y}\right) & z_2 \left(x-\frac{1}{Q_y}\right) & z_3 \left(y-\frac{R}{R_y}\right) \\
    z_1 \left(y-\frac{P}{P_y}\right) & z_2 \left(y-\frac{Q}{Q_y}\right) & z_3 \left(y-\frac{R}{R_y}\right) \\
    z_1                              & z_2                              & z_3 \\
  \end{array}
\right),
$$
where
$$
\mu=\frac{-(P-Q)(Q-R)(R-P)}{[P_y(Q-R)+Q_y(R-P)+R_y(P-Q)]^2}.
$$
and 
$$
\begin{array}{c}
z_1=\frac{P_y(Q-R)}{P_y(Q-R)+Q_y(R-P)+R_y(P-Q)},\ \
z_2=\frac{Q_y(R-P)}{P_y(Q-R)+Q_y(R-P)+R_y(P-Q)},\ \
z_3=\frac{R_y(P-Q)}{P_y(Q-R)+Q_y(R-P)+R_y(P-Q)}.
\end{array}
$$
Thus we have constructed the map
$$
F: V\to SL(3).
$$
 Let us denote the pull-back of the Maurer-Cartan form of $SL(3)$ by $\Omega$:
$$
\Omega:=F^{-1}dF.
$$
This pull-back is called the {\it Darboux derivative} of $F$. The Fundamental Theorem of Calculus for $F$ (see \cite{Sc}) reads as follows. 
\begin{proposition}\label{fundamental} If the Darboux derivatives of two maps $F,\widetilde{F}: V\to SL(3)$ coincide then there is a fixed element $G\in SL(3)$ such that $\widetilde{F}=G\cdot F$
\end{proposition}
Excluding the  case of regular 3-web, one computes
\begin{equation}\label{Omega}
\Omega=\left(
  \begin{array}{ccc}
   \frac{ 1-2a-c}{3}U_3-\frac{1+2a+b}{3}U_2 & bU_2                                     & cU_3 \\
    aU_1                                    & \frac{ 1-2b-a}{3}U_1-\frac{1+2b+c}{3}U_3 & cU_3  \\
    aU_1                                    & bU_2                                     & \frac{ 1-2c-b}{3}U_2-\frac{1+2c+a}{3}U_1\\
  \end{array}
\right),
\end{equation}
where
\begin{equation}\label{Ui}
\begin{array}{cc}
U_1=\frac{[Py(R-Q)+Qy(P-R)+Ry(Q-P)](dy-Pdx)}{(P-R)(P-Q)}, & a=\frac{(P-Q)(P-R)(R-Q)P_{yy}}{[P_y(Q-R)+Q_y(R-P)+R_y(P-Q)]^2}, \\
&\\
U_2=\frac{[Py(R-Q)+Qy(P-R)+Ry(Q-P)](dy-Qdx)}{(Q-P)(Q-R)}, & b=\frac{(Q-R)(Q-P)(P-R)Q_{yy}}{[P_y(Q-R)+Q_y(R-P)+R_y(P-Q)]^2},\\
&\\
U_3=\frac{[Py(R-Q)+Qy(P-R)+Ry(Q-P)](dy-Rdx)}{(R-Q)(R-P)}, & c=\frac{(R-P)(R-Q)(Q-P)R_{yy}}{[P_y(Q-R)+Q_y(R-P)+R_y(P-Q)]^2}.
\end{array}
\end{equation}
Note that the form $U_i$ vanishes on the $i$-th foliation. Moreover, holds
\begin{equation}\label{sum}
U_1+U_2+U_3=0.
\end{equation}

\begin{lemma}
The forms $U_1,U_2,U_3$ and the functions $a,b,c$ are projectively invariant.
\end{lemma}
{\it Proof:}
Identifying locally the projective group  $PGL(3)$ with the linear group $SL(3)$, one observes that the projective  action of some
 $G\in SL(3)$ on the set $V\in \mathbb P^2$ corresponds to the left translation on $SL(3)$. In fact, the normalization (\ref{norm}) implies $F(G\cdot p)=G\cdot F(p)$.
 Since the matrix-valued form $\Omega$ is invariant, so are all its entries.  One checks easily that $U_1=\Omega_{2,2}-\Omega_{3,3}+\Omega_{2,3}-\Omega_{3,2}$, therefore $U_1$ (and, similarly, $U_2$, $U_3$) is invariant. Now the invariance of $U_1$ and of
 $\Omega_{2,1}=aU_1$ implies that $a$ (and, similarly, $b$, $c$) is invariant.
\hfill $\Box$

\section{Linearizability of 3-webs and projective Cartan connections}\label{secstructure}
Following Blaschke, let us define 3 differential operators $\partial_i,$ acting on functions 
$f:V\to \mathbb K$, by
\begin{equation}\label{differentials}
df=f_2U_1-f_1U_2=f_3U_2-f_2U_3=f_1U_3-f_3U_1,\ \ \ f_i=\partial_if.
\end{equation}
\begin{lemma}\label{structureEQlemma}
The forms $U_1,U_2,U_3$ and the functions $a,b,c$ satisfy the following structure equations:
\begin{equation}\label{structureEQ}
\begin{array}{lll}
dU_1=(c-b)U_2\wedge U_3, &  dU_2=(a-c)U_3\wedge U_1, & dU_3=(b-a)U_1\wedge U_2,\\
\\
a_1=a[1+2(b-c)],  & b_2=b[1+2(c-a)],  & c_3=c[1+2(a-b)].
\end{array}
\end{equation}
These equations are equivalent to one matrix equation \begin{equation}\label{MC}
d\Omega +\Omega\wedge \Omega=0.
\end{equation}
\end{lemma}
{\it Proof:}
The matrix-valued form $\Omega$, being the pull-back of the Maurer-Cartan form, verifies the structure equation (\ref{MC}).
Substituting the expression (\ref{Omega}) for $\Omega$, taking into account the identity (\ref{sum}) and the definition (\ref{differentials}), one sees that equations (\ref{structureEQ}) are equivalent to  one matrix equation (\ref{MC}).
\hfill $\Box$\\

\noindent{\bf Remark 1.} Observe that the system (\ref{structureEQ}) is symmetric with respect to an action of the permutation group $S_3$. For example, the transposition $(1,2)$ acts as follows:
$$
U_1\mapsto -U_2, \ \ U_2\mapsto -U_1,\ \ U_3\mapsto -U_3,\ \ a\mapsto -b,\ \ b\mapsto -a,\ \ c\mapsto -c.\ \
$$

Suppose that a planar 3-web is described by three 1-forms $U_i$, each vanishing on leaves of its "own" foliation $\mathcal{F}_i$, and the forms are normalized to satisfy (\ref{sum}). Note, that this normalization is determined up to rescaling by a non-vanishing factor $U_i \to \frac{1}{N}U_i$.
\begin{lemma}\label{linLM}
If the forms $U_1,U_2,U_3=-(U_1+U_2)$   satisfy equations
(\ref{structureEQ}), then the web is linearizable.
\end{lemma}
{\it Proof:} Let us fix some point $p_0$, define the matrix-valued differential form  $\Omega$ by (\ref{Omega}), and consider the matrix Pfaff equation
\begin{equation}\label{automorph}
dF=F\cdot\Omega.
\end{equation}
This equation is integrable due to Lemma \ref{structureEQlemma}. Therefore for any $G_0\in SL(3)$ there is a unique solution $F(p)$ with the initial condition $F(p_0)=G_0$. Let $\xi_i(p)$ be the columns of this solution  $F(p)=(\xi_1(p),\xi_2(p),\xi_3(p))$ and $\zeta(p):=\xi_1(p)+\xi_2(p)+\xi_3(p)$. We claim that the map
$p\mapsto \eta(p):=[\zeta(p)]\in \mathbb P^2$ linearizes the web.  In fact, equations (\ref{differentials}) give $U_1(\partial_1)=0$, $U_2(\partial_1)=-1$, $U_3(\partial_1)=1.$ One computes $\partial_1(\zeta)=\xi_1-\frac{1}{3}(1+2(b-c))\zeta.$ Therefore the tangent line to the leaf of $\mathcal{F}_1$ through $\eta(p)$  is spanned by $[\zeta(p)]$ and $[\xi_1(p)]$. Since $\partial_1(\xi_1)=\frac{1}{3}(2+b-c)\xi_1$ this tangent line is stable along the leaf. Hence the constructed map rectifies the foliation $\mathcal{F}_1$. Due to the symmetry of equations (\ref{structureEQ}) our map linearizes also $\mathcal{F}_2$ and $\mathcal{F}_3$.
\hfill $\Box$\\

\noindent{\bf Remark 2.} The system of uncoupled Euler equations (\ref{Euler}) is invariant with respect to the action of the projective group $PGL(3)$ in the $xy$-plane of independent variables, prolonged on the inclinations $P,Q,R$ to preserve the distributions $dy-Pdx=dy-Qdx=dy-Rdx=0$. Given $\Omega$, the group $PGL(3)$ acts  transitively on the space of solutions to
the matrix equation  (\ref{automorph}).   Thus, the pair of equations (\ref{automorph},\ref{structureEQ}) is the so-called {\it group splitting} of (\ref{Euler}) into the  {\it automorph} system (\ref{automorph}) and the {\it resolving} system (\ref{structureEQ}) (see \cite{Vi,Og}). Any $\Omega $, defined by a solution to (\ref{structureEQ}), labels some $PGL(3)$-orbit of a solution to (\ref{Euler}). The orbit itself is the space of solutions to (\ref{automorph}).

\begin{lemma}\label{equalInv}
Suppose that a diffeomorphism $\varphi:V\to \widetilde{V}$ maps a linear 3-web $\mathcal{L}_3$ into a linear 3-web $\widetilde{\mathcal{L}}_3$. If the Darboux derivatives of their corresponding maps $F, \widetilde{F}$ verify $\varphi^*(\widetilde{\Omega})=\Omega$ then the webs are projectively equivalent and there is $G\in PGL(3)$ such that $\varphi|_V=G$.  
\end{lemma}
{\it Proof:} Pulling back $d\widetilde{F}$ by $\varphi$ we get $d(\widetilde{F}\circ\varphi)=\varphi^*(d\widetilde{F})=\varphi^*(\widetilde{F}\cdot \widetilde{\Omega})=\widetilde{F}\circ\varphi\cdot \varphi^*(\widetilde{\Omega})=\widetilde{F}\circ\varphi\cdot \Omega$. Therefore the Darboux derivatives of $F$ and $\widetilde{F}\circ\varphi$ coincide and by Proposition \ref{fundamental} holds true $\widetilde{F}\circ\varphi=G\cdot F$ for some $G\in PGL(3)$. Invoking the construction of maps $F$ and $\widetilde{F}$, we conclude that the webs are projectively equivalent.
\hfill $\Box$\\

\noindent The Chern connection form $\gamma$, defined by $dU_i=\gamma \wedge U_i$, in our normalization reads as
$$
\gamma =aU_1+bU_2+cU_3.
$$
Using (\ref{structureEQ}), we get the Blaschke curvature form:
\begin{equation}\label{curvature}
d\gamma =(a+b+c)U_1\wedge U_2.
\end{equation}

Let a planar non-hexagonal 3-web be determined on some open set $V\in \mathbb C^2$ by three 1-forms $\omega_1,\omega_2,\omega_3=-(\omega_1+\omega_2)$. Then one can renormalize these forms (see \cite{BB}) so that the Blaschke curvature is $\omega_1\wedge
\omega_2$ . In this normalization 
\begin{equation}\label{cur1}
d\omega_1=\alpha\omega_1\wedge
\omega_2, \ \ d\omega_2=\beta\omega_1\wedge
\omega_2, \ \ \ \beta_1-\alpha_2=1,
\end{equation}
where the notation (\ref{differentials}) is used for differentiation with respect to $\omega_i$. Choose four functions $a,b,c,N$ and construct the $\mathfrak{sl}(3)$-valued form (\ref{Omega}) with $U_i=\frac{\omega_i}{N}$. Conceptually, the form $\Omega$ defines a {\it projective Cartan connection} by specifying the {\it Cartan gauge} $(V,\Omega)$ (see \cite{Sc}, page 174).
 This connection has the curvature form 
\begin{equation}\label{projcur}
K=d\Omega +\Omega\wedge \Omega.
\end{equation}
Let us try to adjust the choice of $a,b,c,N$ to kill the curvature $K$. 
Analysis of 9 scalar equations $K=0$ quickly gives the following expressions for $a,b,c$: 
\begin{equation}\label{abc}
\begin{array}{l}
a=\frac{N^2}{3}+\left(\frac{2\beta}{3}+\frac{\alpha}{3}\right)N-\frac{N_1}{3}-\frac{2N_2}{3}\\
\\
b=\frac{N^2}{3}-\left(\frac{\beta}{3}+\frac{2\alpha}{3}\right)N+\frac{2N_1}{3}+\frac{N_2}{3}\\                                  
\\
c=\frac{N^2}{3}+\left(\frac{\alpha}{3}-\frac{\beta}{3}\right)N-\frac{N_1}{3}+\frac{N_2}{3}
\end{array}
\end{equation}
Substituting these expressions again into (\ref{projcur}), one computes 
\begin{equation}\label{curform}
K=\left(
  \begin{array}{ccc}
   K_{11},  & K_{22}, &  -K_{11} - K_{22}\\
   K_{11},  & K_{22}, &  -K_{11} - K_{22}\\
   K_{11},  & K_{22}, &  -K_{11} - K_{22}\\
  \end{array}
\right)\omega_1\wedge
\omega_2,
\end{equation}
where
$$
\begin{array}{l}
K_{11}=  \frac{2N_1^2+4N_1N_2+N_1+2N_2}{3N^2}-\frac{N_{11}+N_{12}+N_{21}}{3N}-\frac{(\alpha+\beta)N_1+\alpha N_2+\alpha+2\beta}{N}+
 \frac{2\alpha N}{3}+\frac{2\alpha^2+4\alpha\beta+\alpha_1+\alpha_2+\beta_1}{3},\\
 \\
K_{22}= \frac{N_{22}+N_{12}+N_{21}}{3N} -\frac{2N_2^2+4N_1N_2+2N_1+N_2}{3N^2}+\frac{\beta N_1+(\alpha+\beta)N_2+2\alpha+\beta}{N}+
 \frac{2\beta N}{3}-\frac{2\beta^2+4\alpha\beta+\beta_1+\beta_2+\alpha_2}{3}.

\end{array}
$$
The group $\rm Aff(2)$ of affine transformations of the 2-dimensional plane can be realized as the $SL(3)$ stabilizer of $[1:1:1]\in \mathbb P^2$, the corresponding sub-algebra $ \mathfrak{aff}(2)\subset \mathfrak{sl}(3)$ annihilating the vector $(1,1,1)^T$    
\begin{theorem}\label{linearization}
Let the forms $w_i$ of a planar non-hexagonal 3-web be normalized 
as in (\ref{cur1}) and $N$ be a non-vanishing function. Then the form $\Omega_N$ constructed as in (\ref{Omega}) with $U_i=\frac{\omega_i}{N}$ and $a,b,c$ as in (\ref{abc}) defines a torsion-free projective Cartan connection with the model geometry $(\mathfrak{sl}(3),\mathfrak{aff}(2))$,  the web leaves being its geodesics. The web is linearizable  if and only if there is  
$N$ for which this connection is flat: $K_{11}=K_{22}=0.$ 
\end{theorem}
{\it Proof:} One checks that the linear map  
$$
 T_pV \stackrel{\Omega}{\longrightarrow}  \mathfrak{sl}(3) \to \mathfrak{sl}(3)/ \mathfrak{aff}(2)
$$
is an isomorphism. Therefore $\Omega_N$ defines a Cartan gauge with the model geometry $(\mathfrak{sl}(3),\mathfrak{aff}(2))$. By formula (\ref{curform}), the curvature $K$ takes values in $\mathfrak{aff}(2)$. Therefore the corresponding  projective Cartan connection is torsion-free (for definitions and details see \cite{Sc}).

If $\partial=\frac{d}{ds}$ is the differentiation along some parametrized curve $C$ in $V$ then its development 
$s\mapsto F(s)=(\xi_1(s),\xi_2(s),\xi_3(s))\in SL(3)$ is the solution to $\partial F=F\cdot\Omega_N$ with $F(0)=e$. The curve $c$ is a geodesic if $\zeta(s)=F(s)\cdot(1,1,1)^T\in \mathbb K^3$ represents a line in $\mathbb P^2$. If $C$ is a leave of, say, the first foliation then choosing $\partial=\partial_1$ one has $\partial_1(\zeta)=\xi_1(s)-\frac{1}{3}(1+2(b-c))\zeta(s)$. Therefore the tangent to $[\zeta(s)]$ is spanned by $[\zeta(s)]$ and $[\xi_1(s)]$. Since $\partial_1(\xi_1)=\frac{1}{3}(2+b-c)\xi_1$ this tangent is stable along the leaf (compare with the calculations in the proof of Lemma \ref{linLM}) and the leaf is geodesic. 

\noindent Finally, by Lemmas \ref{structureEQlemma} and \ref{linLM},   the web is linearizable if and only if $K=0$. 
\hfill $\Box$
\smallskip

\noindent{\bf Remark 3.} Given a planar non-hexagonal 3-web, the system $K_{11}=K_{22}=0$ is overdetermined: we have two second order partial differential equations for one unknown function $N$. A compatibility analysis quickly gives all second order derivatives of $N$ in terms of the web invariants $\alpha,\beta$ and their derivatives up to the second order. Then the conditions $d(dN_1)=d(dN_2)=0$ give two equations, quadratic in $N_1,N_2$. These two compatibility conditions and the expressions for $N_{ij}$ are surprisingly short: they would easily fit in half page. (We do not give them as they will not be used.) However, the further analysis is possible only with some computer algebra software.
Therefore the idea to write down  the compatibility conditions explicitly seems rather unpromising.
As follows from the known results (see the discussion in Introduction), this analysis would give a polynomial equation for $N$ of degree at most 15. Thus, a planar non-hexagonal 3-web        
carries a finite number of "natural" projective Cartan connections, determined implicitly by the polynomial. In contrast, a planar hexagonal 3-web carries a one-parameter family of flat projective Cartan connections.

\section{Signature sets}
Let a linear planar 3-web be defined on some open set $V\subset \mathbb P^2$. For any non-regular web, formulas (\ref{Ui},\ref{differentials}) define a map $\sigma:V \to \mathbb K^9$, $p\mapsto\ (a,b,c,a_2,b_3,c_1,a_{22},b_{33},c_{11})$, where $a_{22}=\partial_2a_2$, $b_{33}=\partial_3b_3$ and $c_{11}=\partial_1c_1$.
\begin{definition}
The signature set $\mathcal{S}_{\scriptscriptstyle \mathcal{L}_3}$ of a linear planar non-regular 3-web $\mathcal{L}_3$ is the image of $V$ under the above defined map, i.e. $\mathcal{S}_{\scriptscriptstyle \mathcal{L}_3}:=\sigma(V)$.
 \end{definition}
One expects that generically the map $\sigma$ parametrizes some (possibly singular) surface in $\mathbb K^9$. However, for some webs the signature set degenerates to a (possibly singular) curve or even to a point. For example, a linear web is formed by 3 pencils of straight lines if and only if $a=b=c=0$, which is equivalent to $P_{yy}\equiv Q_{yy}\equiv R_{yy}\equiv 0$ for the corresponding solution to (\ref{Euler}) (see Remark after Lemma \ref{aligned}).
\begin{theorem}\label{0dim}
If the signature set $\mathcal{S}_{\scriptscriptstyle \mathcal{L}_3}$ of a linear planar non-regular 3-web consists of one point, then either $\mathcal{S}_{\scriptscriptstyle \mathcal{L}_3}=\{(0,0,0,0,0,0,0,0,0)\}$ and the web is formed by 3 pencils of straight lines; or one can enumerate the web foliations so that $\mathcal{S}_{\scriptscriptstyle \mathcal{L}_3}=\{(\frac{1}{2},-\frac{1}{2},0,0,0,0,0,0,0)\}$ and the web is formed by tangents to a conic and by a pencil of lines centered on this conic.
\end{theorem}
{\it Proof:} Since $a,b,c$ are constant, it is immediate that $a_2=b_3=c_1=0$. The second line of equations (\ref{structureEQ}) implies $a[1+2(b-c)]=b[1+2(c-a)] =c[1+2(a-b)]=0$, which gives $\mathcal{S}_{\scriptscriptstyle \mathcal{L}_3}$ as announced above. The case when $\mathcal{S}_{\scriptscriptstyle \mathcal{L}_3}$ sits in the origin was considered above. Computing the invariants $a,b,c$ for the web formed by tangents to a parabola and by lines parallel to its axis, one obtains  $a=\frac{1}{2},\ b=-\frac{1}{2},\ c=0.$ Now all webs with this signature set have the same structure equations for $U_i$  (see the first line of (\ref{structureEQ})), and the theorem follows from Lemma \ref{equalInv}.
\hfill $\Box$\\

The degeneration of the signature set is explained by projective symmetries of the web.  
\begin{definition}
An infinitesimal symmetry of a d-web is a vector field whose local flow preserves the web.
 \end{definition}
In fact, both types of webs, described by Theorem \ref{0dim}, possess 2-dimensional projective symmetry algebras (see \cite{Agr} for the classifications of linear 3-webs admitting infinitesimal symmetries).

Now let us described the webs whose signature set is one-dimensional. We will need the following Lemma. 

\begin{lemma}\label{formwithSym}
A differential form $\omega=p(u,v)du+q(u,v)dv$ is invariant along the local flow of a vector field $X=\lambda(u,v)\partial_v$ if and only if holds
\begin{equation}\label{formwithSymEQ}
\partial _v\left(\frac{p_v}{q}\right)=\partial _u\left(\frac{q_v}{q}\right).
\end{equation}
\end{lemma}
{\it Proof:} The form is invariant if and only if its Lie derivative $\mathcal{L}_X(\omega)=(\lambda p_v+q\lambda_u)du+(\lambda q_v+q\lambda_v)dv$ vanishes. 
Therefore $q\lambda_u=-\lambda p_v, \ \ q\lambda_v=- \lambda q_v.$ Hence (\ref{formwithSymEQ}).
\hfill $\Box$\\

\begin{theorem}\label{1dim}
The signature set of a linear planar 3-web is one-dimensional if and only if the web admits a one-dimensional symmetry group of projective transformations.
\end{theorem}
{\it Proof:} If the web is symmetric with respect to one-dimensional subgroup of projective transformations then the projective invariants 
$(a,b,c,a_2,b_3,c_1,a_{22},b_{33},c_{11})$ are constant along the orbits and therefore $\dim \mathcal{S}_{\scriptscriptstyle \mathcal{L}_3}\le 1$. The signature set cannot degenerate to a point since such webs possess 2-dimensional projective symmetries. 

Now suppose that $\dim \mathcal{S}_{\scriptscriptstyle \mathcal{L}_3}= 1$  and $\mathcal{S}_{\scriptscriptstyle \mathcal{L}_3}$ is parametrized by one parameter $u:V\to \mathcal{S}_{\scriptscriptstyle \mathcal{L}_3}$. Choose a function $v$ so that $(u,v)$ are local coordinates, and,  for $U_1=pdu+qdv, \ U_2=mdu+ndv$ with some functions $p,q,m,n$,  the coefficient $q$ do not vanish. From
$da=a'(u)du=a_2U_1-a_1U_2$ we get $\frac{a_1}{a'}=\frac{q}{\Delta}$, where $\Delta=pn-mq$. Due to the equation $a_1=a[1+2(b-c)]$, the function $a_1$ also depends only on $u$. Therefore holds  $\partial _v\left(\frac{q}{\Delta}\right)=0.$ Similarly, from analysis of $db$ one gets $\frac{b_2}{b'}=\frac{n}{\Delta}$ and  $\partial_v\left(\frac{n}{\Delta}\right)=0.$ With $s=\ln \Delta$ we derive $q_v=qs_v,\ n_v=ns_v$ and $q_{uv}=q_us_v+qs_{uv}$.
The structure equation $dU_1=(c-b)U_1\wedge U_2$ implies $c-b=\frac{q_u-p_v}{\Delta}$ and therefore $\partial_v\left( \frac{q_u-p_v}{\Delta}\right)=0$. Differentiating and taking into account the above found $q_{uv}$, we calculate
$p_{vv}=s_vp_v+qs_{uv}$ and (\ref{formwithSymEQ}) follows. Due to Lemma \ref{formwithSym}, there is 
a vector field $X=\lambda(u,v)\partial_v$ with $\lambda$ defined up to a constant by $q\lambda_u=-\lambda p_v, \ \ q\lambda_v=- \lambda q_v$, whose local flow leaves $U_1$ invariant. Let us show that also $\mathcal{L}_X(U_2)=0$. Applying the Lie derivative $\mathcal{L}_X$ to $dU_1=(c-b)U_1\wedge U_2$ we get $\mathcal{L}_X(U_1\wedge U_2)=0$ and therefore $\mathcal{L}_X(U_2)=\nu  U_1$. In coordinates one has 
$\mathcal{L}_X(U_2)\wedge du=(\lambda n_v+n\lambda_v)dv\wedge du=(\lambda ns_v-n\frac{\lambda q_v}{q})dv\wedge du=(\lambda ns_v-\lambda ns_v)dv\wedge du=0=\nu U_1\wedge du=q \nu dv\wedge du$. As $q\ne0$ the last equality implies $\nu=0$ and $\mathcal{L}_X(U_2)=0$.
Finally $\mathcal{L}_X(U_3)=\mathcal{L}_X(-U_1-U_2)=0$ and $X$ is an infinitesimal symmetry of the web. Any transformation $exp(tX)$ of the local flow leaves invariant $U_i$ and $a,b,c$. Thus $\Omega$ is invariant and $exp(tX)$ is projective by Lemma \ref{equalInv}. 
\hfill $\Box$

\begin{lemma}\label{propformsforcurve}
If two linear planar 3-webs have the same one-dimensional signature set then one can choose the local coordinates so that the forms $U_i$ of the webs coincide in the chosen coordinates.   
\end{lemma}
{\it Proof:} Suppose that our two webs $\mathcal{L}_3$ and $\tilde{\mathcal{L}}_3$ are defined on open sets $V$ and $\tilde{V}$ and have the same one-dimensional signature set $\mathcal{S}_{\scriptscriptstyle \mathcal{L}_3}$.  Let us parametrize the curve $\mathcal{S}_{\scriptscriptstyle \mathcal{L}_3}$ by some parameter. Then this parameter pulls back to $V$ and $\tilde{V}$ and define there functions $u$ and $\tilde{u}$. Let us choose the functions $v$, $\tilde{v}$ on the sets $V,\tilde{V}$ so that the infinitesimal symmetries of the webs assume the forms $\partial_v$ and $\partial_{\tilde{v}}$ respectively. The pairs $(u,v)$ and $(\tilde{u},\tilde{v})$ give local coordinate systems.
For each of the webs, at least two of the three forms $U_i$ and at least two of the three forms $\tilde{U}_i$ have non-vanishing coefficients of $dv$ and $d\tilde{v}$ respectively. Therefore at least for one index $i\in\{1,2,3\}$, say $i=1$, the coefficients of $dv$ and $d\tilde{v}$  in $U_1$ and $\tilde{U}_1$  do not vanish.

For the web $\mathcal{L}_3$, the invariants $a,b,c$ and the coefficients of the forms $U_i$ do not depend  on $v$. Therefore $U_1=p(u)du+q(u)dv$ with  $q(u)\ne 0$. One can change the second coordinate by $v\to v+\chi(u)$ to kill the coefficient $p(u)$. Let $U_2=m(u)du+n(u)dv$, then $U_1\wedge U_2=-mqdu\wedge dv$.  With $da\wedge U_1=a_1U_1\wedge U_2$, $db\wedge U_1=b_2U_1\wedge U_2$ and the equation for $dU_1$ we obtain: 
\begin{equation}\label{formsforcurve}
\frac{da}{du}=-ma_1,\ \ \ \ \  n\frac{db}{du}=-mqb_2,\ \ \ \ \  \frac{dq}{du}=(b-c)mq.
\end{equation}
Note that $a_1=a[1+2(b-c)]$ and $b_2=b[1+2(c-a)]$ are also the functions only of $u$. Therefore $m$ is completely determined by the signature set,  $q$  is defined up to a constant factor, and, finally, this factor fixes $n$. 

Deriving the counterpart of (\ref{formsforcurve}) for $\tilde{\mathcal{L}}_3$ and rescaling, if necessary, the coordinate $\tilde{v}$, we make the forms $U_i$ coincide with  corresponding $\tilde{U}_i$ in the obtained coordinates.    \hfill $\Box$\\

\noindent{\bf Remark 1.} A one-dimensional infinitesimal symmetry $X$ of a linear non-hexagonal 3-web  is projective.
In fact, if  $exp(tX)$ is the local flow of the symmetry then each $t$ gives a map $exp(tX)$ respecting the linearity of the web. 
For a non-hexagonal linear 3-web, there are only a finite number of such maps that are projectively non-equivalent (see \cite{Bk}), thus  $exp(tX) \in PGL(3)$. A complete classification of linear non-hexagonal 3-webs with one infinitesimal symmetry was obtained in \cite{Agr}. Moreover, there was presented a classification of linear hexagonal 3-webs with infinitesimal {\it projective} symmetries. 

\begin{theorem}\label{signature}
Suppose that the signature sets of two linear planar non-regular 3-webs coincide in a neighborhood of a non-singular point. Then the web germs are projectively equivalent.
\end{theorem}
{\it Proof:}  If the signature sets are points,  and these points coincide, then the webs are projectively equivalent due to Theorem \ref{0dim}.
If the signature set is a curve then the claim follows from  Lemmas \ref{propformsforcurve} and \ref{equalInv}. 

Finally, if the signature set is 2-dimensional then two of the invariants $a,b,c,a_2,b_3,c_1,$ $a_{22},b_{33},c_{11}$ can be chosen as local coordinates. Now the other 7 invariants and all their derivatives are functions of the chosen two.  Therefore the forms $U_1,U_2$ are uniquely defined by the signature set. (For example, if $da\wedge db\neq 0$, then $c,a_2,b_3,c_1,a_{22},b_{33},c_{11}$ are functions of $a,b$ and by
$da=a_2U_1-a[1+2(b-c)]U_2$, $db=b[1+2(c-a)]U_1 +(b[1+2(c-a)]+b_3)U_2$ the forms $U_1,U_2$ are uniquely defined. We have used $b_1=-b_2-b_3$.)  
Now the form $\Omega$ is the same for our two webs, and the webs are projectively equivalent by Lemma \ref{equalInv}.\hfill $\Box$\\

\noindent{\bf Remark 2.} While the condition $a=b=c\equiv0$ distinguishes 3-webs of 3 pencils of lines, a simple relation  $a+b\equiv 0$ (or $b+c\equiv 0$, or $c+a\equiv 0$) characterizes 3-webs, whose 2 foliations are formed by tangents to one and the same conic. In fact, the relation  $a+b\equiv 0$ is equivalent to $P_{yy}+Q_{yy}\equiv 0$.  Let us replace our web $\mathcal{L}_3$ by a 3-web $\widetilde{\mathcal{L}}_3$, whose 2 foliations are the same as described by $P(x,y),Q(x,y)$ and the third one is some pencil of lines. For the invariants $\tilde{a},\tilde{b},\tilde{c}$ of $\widetilde{\mathcal{L}}_3$ holds true $\tilde{a}+\tilde{b}\equiv0$ and $\tilde{c}\equiv 0$. Therefore $\tilde{a}+\tilde{b}+\tilde{c}\equiv 0$ and the web $\widetilde{\mathcal{L}}_3$ is hexagonal. By the classical result of Graf and Sauer \cite{GSg}, the lines of $\widetilde{\mathcal{L}}_3$ are tangent to a curve of 3d class (i.e. to the dual of some cubic). By construction of $\widetilde{\mathcal{L}}_3$, this curve degenerates to a point (namely, the pencil center) and a conic. 
Conversely, one  verifies easily that $P_{yy}+Q_{yy}\equiv 0$ is true for the inclinations $P(x,y),Q(x,y)$ of two tangents to a conic, passing through a point $(x,y)$.\\

\noindent{\bf Remark 3.} For non-symmetric webs with non-constant $a,b,c$, one does not need 9-dimensional space to define the signature set and may reduce the number of invariants to three, namely one can choose $a,b,c$. Their derivatives are needed for webs whose dual focal curves are lines or belong to the same conic. For example,  projective orbits of linear 3-webs with $a=b=0$ (i.e. whose 2 foliations are pencils of lines) need invariants $c_1$ and $c_{11}$ to be separated.

\section{Polymorphic 3-webs} 
In this section we use the obtained invariant description to deduce some properties of polymorphic planar 3-webs. 
First of all, to control the hexagonality of the web, which is equivalent to $k:=a+b+c\equiv0$, we rewrite the structure equations in terms of $(a,b,k)$ and their derivatives with respect to $\partial_1$ and $\partial_2$: 
\begin{equation}\label{structureEQk}
\begin{array}{c}
dU_1=(k-a-2b)U_1\wedge U_2, \ \ \ \ \ \ dU_2=(2a+b-k)U_1\wedge U_2,\\
\\
a_1=a[1+2(a+2b-k)], \ \ \ \ \ \ \ b_2=b[1+2(k-2a-b)],  \\
\\
k_1+k_2=a_2+b_1+2(a+b)+4(a^2-b^2)+4k(b-a)-k.
\end{array}
\end{equation}
To write the equation for $k$ we have used the identity $\partial_1+\partial_2+\partial_3=0$.\\

\noindent{\bf Remark 1.} It follows from equations (\ref{structureEQk}) that $k=const$ (i.e. $k_1=k_2=0$) implies $k=0$. One can check this as follows. Introducing a new parameter $m$ by $2m=b_1-a_2$ we express $b_1$ and $a_2$ via $m$ from the last equation of (\ref{structureEQk}). Now the compatibility conditions $d(da)=d(db)=0$ give $m_1$ and $m_2$. Then the equation $d(dm)=0$ gives $m$, provided that $k\ne0$. Differentiating $m$ and comparing $m_1,m_2$ with the expressions obtained earlier, we get two independent polynomial equations, involving $a,b,k$. They imply that $a,b$ are also constant. Therefore $a_1=b_2=c_3=0$. From (\ref{structureEQ}) we have $a_1+b_2+c_3=a+b+c$, hence $k=0$.\\

Suppose that a linear web $\mathcal{L}_3$ is polymorphic, i.e. there is a non-projective map $\varphi: V\to \mathbb P^2$ respecting the linearity. Let $\widetilde{U}_i$ be the invariant forms (\ref{Ui}) of the transformed web $\varphi(\mathcal{L}_3)$. Thus, for the pull-backs we have $\varphi^{*}\widetilde{U}_i=(1+f)U_i$, where $f\not\equiv 0$ since the map is not projective (see Lemma \ref{equalInv}). Abusing notation, we can think of the pull-backs $\varphi^{*}\widetilde{U}_i$ as of the re-normalization $\widetilde{U}_i=(1+f)U_i$, $i=1,2,3$ of the forms $U_i$. Since the transformed web $\varphi(\mathcal{L}_3)$ is linear, the re-scaled forms $\widetilde{U}_i$ also satisfy equations (\ref{structureEQ}).
 The invariants of the alternative linear form  $\varphi(\mathcal{L}_3)$ of the  web $\mathcal{L}_3$ are as follows:
\begin{equation}\label{transformedKab}
\begin{array}{lll}
\widetilde{a}=\frac{f_1+2f_2+(3a-k)f+3a}{3(f+1)^2}, & \widetilde{b}=\frac{-2f_1-f_2+(3b-k)f+3b}{3(f+1)^2}, & \widetilde{k}=\frac{k}{(f+1)^2},
\end{array}
\end{equation}
where the sub-indices denote, as before, the derivations by $\partial_i$. 
\begin{lemma}\label{polymorph}
A linear web with the structure equations (\ref{structureEQk}) is polymorphic if and only if there is a non-vanishing solution $f$ of the following system: 
\begin{equation}\label{polymorphEQ}
\begin{array}{l}
f_{11}+f_{12}+f_{21}=f(f_1+2f_2)+[1+3(b-a)]f_1+[2+3(a+2b-k)]f_2+\\
\\
\ \ \ \ \ \ (3a-k)f^2+[2k^2-2k(a+2b)+k_1-k+3a]f,\\
\\
f_{22}+f_{12}+f_{21}=f(2f_1+f_2)+[2+3(k-2a-b)]f_1+[1+3(b-a)]f_2+\\
\\
\ \ \ \ \ \ (k-3b)f^2+[2k^2-2k(2a+b)-k_2+k-3b]f.\\

\end{array}
\end{equation}
\end{lemma}
{\it Proof:} For the differentiations with respect to the rescaled forms $\widetilde{U}_i$ one has $\widetilde{\partial}_i=\frac{1}{1+f} \partial_i$. The invariants
$\widetilde{a},\widetilde{b},\widetilde{k}$ satisfy 
$$
\widetilde{\partial}_1(\widetilde{a})=\widetilde{a}[1+2(\widetilde{a}+2\widetilde{b}-\widetilde{k})], \ \ \ \ \ \ \ \widetilde{\partial}_2(\widetilde{b})=\widetilde{b}[1+2(\widetilde{k}-2\widetilde{a}-\widetilde{b})]. 
$$
These two equations are equivalent to (\ref{polymorphEQ}). The equation for $\widetilde{k}$, corresponding to the last equation of (\ref{structureEQk}), follows from (\ref{polymorphEQ}).

If there is a non-vanishing solution $f$ to (\ref{polymorphEQ}) then the 3-web with the invariants defined by (\ref{transformedKab}) and the rescaled forms $\widetilde{U}_i$ admits a linearization by Lemma \ref{linLM}.   The condition $f\not\equiv 0$ ensures that this linearization is not projectively equivalent to the identity.
\hfill $\Box$\\

\noindent{\bf Remark 2.} One may be tempted to search for particularly simple solutions to (\ref{polymorphEQ}), for instance, such that $f=const\ne 0$. Unfortunately, this Ansatz does not work.  Since all derivatives of $f$ vanish, equations (\ref{polymorphEQ}) give $k_1,k_2$. Then the equation $d(dk)=0$ and the last equation of (\ref{structureEQk}) determine $a_2,b_1$. Thus, all first derivatives of $a,b,k$ are expressed in terms of $a,b,k$ and of the constant $f$. The equations  $d(da)=d(db)=0$ and their derivatives give 6 polynomial equations for $a,b,k$, incompatible with $k\ne 0$.\\

Due to the nonlinearity of the weakly overdetermined system (\ref{polymorphEQ}), its compatibility analysis  is impossible without help of symbolic computation software. The usual approach was to exclude the functions, defining the linearizing map (i.e. $f$ and its derivatives). This leads very quickly to very involved expressions, unmanageable even by computer algebra. We find more promising to unite equations (\ref{structureEQk}) and (\ref{polymorphEQ}), and consider them as a system characterizing polymorphic 3-webs.
\begin{theorem}\label{analyticdim}
Any polymorphic linear 3-web is analytic. The space of projective moduli of polymorphic linear 3-webs is at most 8-dimensional.
\end{theorem}
{\it Proof:} By classical result of Graf and Sauer, linear hexagonal 3-webs are analytic, and their projective moduli space coincides with that of planar cubics, i.e. it is one-dimensional. Thus, it is enough to consider non-hexagonal webs and set $k\ne 0$. We give here the sketch of the proof and present the details of the computation scheme in the Appendix.

Let us introduce invariant parameters $L$ and $h$ by $2L=f_{12}+f_{21}$, $2hk=k_2-k_1$. Then the differentials of the following 12 invariants $a,a_2,a_{22},b,b_1,b_{11},k,h,f,f_1,f_2,L$ can be expressed  in the form $dI_a=F_a^1U_1+F_a^2U_2$, where $I_a$ are these invariants and the coefficients $F_a^i$ are rational functions of these 12 invariants. The compatibility conditions $d(dI_a)=0$ are not satisfied identically. They give 2 polynomial equations for the invariants $I_a$. If the signature set is one-dimensional then the web has one-dimensional symmetry by Theorem \ref{1dim}. Therefore it is hexagonal, since the Gronwall conjecture is true for such webs (see \cite{Agr}). Thus, the signature set is 2-dimensional and we can choose two of the above 12 invariants, say $I_{\alpha}, I_{\beta}$, as local coordinates and express $U_1,U_2$ in terms of $dI_{\alpha}, dI_{\beta}$. Now the differentials of the left 10 invariants can be written via $dI_{\alpha}, dI_{\beta}$. We obtain a polynomial exterior differential system with 2 constraints. Hence, if there is a (local) solution to it  then this solution is analytic and depends  on at most 8 constants.
\hfill $\Box$\\

\noindent{\bf Remark 3.} The derivatives of the two constraints, obtained by calculation of compatibility conditions in the proof of Theorem \ref{analyticdim}, give 4 more polynomial constraints. There is an evidence that at least 5 of the 6 constraints are independent. Therefore the projective moduli space of polymorphic 3-webs is at most 5-dimensional. But we are 
unable to check the independence with the computational resources available.

\section{Projective invariants and Gronwall's conjecture}

The developed theory permits one to resolve the Gronwall conjecture algorithmically, provided that  sufficiently powerful computational capacity is available. 

In the proof of Theorem \ref{analyticdim}, we explain how to obtain two  constraints $\Phi=\Psi=0$, where $\Phi$ and $\Psi$ are polynomials in 12 invariants $a,a_2,a_{22},b,b_1,b_{11},k,h,f,f_1,f_2,L$.
Let us introduce two new variables $S,T$ and consider the ascending chain of ideals $J_0\subset J_1\subset J_2 \subset ... \subset J_k\subset J_{k+1}...\subset... $, where  $J_0=\langle \Phi,\Psi,kS-1,fT-1 \rangle$ and $J_{k+1}$ is obtained from $J_k$ as follows:
differentiate all the generators of $J_k$ but $kS-1,fT-1$ with respect to $U_1,U_2$, clear the denominators and add the obtained polynomials to the generators  of $J_k$.  

Then if the descending sequence of natural numbers $dim(J_0)\ge dim(J_1)\ge dim(J_2)\ge... \ge dim(J_k)\ge dim(J_{k+1})\ge ...$ stabilizes for some $l$ grater then 1: $dim(J_{l})=dim(J_{l+1})\ge 2$ then the conjecture is false. Indeed, we choose an irreducible component of maximal dimension $d\ge 2$ of the intersection of affine algebraic varieties  $X(J_{l})\cap X(J_{l+1})$ and project it along $ST$-plane, thus obtaining an affine algebraic variety in $\mathbb C^{12}$ equipped with a consistent polynomial exterior differential system, defining a polymorphic non-hexagonal 3-web.  

If  $dim(J_{l})$ jumps at some step for a value less then 2 then the conjecture is true. In fact, for the non-empty algebraic set $X(J_{l})$, the signature set $\mathcal{S}_{\scriptscriptstyle \mathcal{L}_3}$ would be either  a curve, and the  web would have a one-dimensional symmetry, or a point, and the web would be hexagonal. By the main result of \cite{Agr}, the conjecture is true for webs with infinitesimal symmetries.

In this section we show how this scheme works for 3-webs whose 2 foliations are 2 pencils of lines. 
To perform the calculations described in this section, one needs  a symbolic computation software. The author used Maple 18 installed on a computer with 16GB of memory. 

\subsection{3-webs with two pencils of lines} 
As we have shown, for such webs holds $a\equiv b \equiv 0$, hence $k=c$. First we rewrite equations (\ref{structureEQk},\ref{polymorphEQ}) in a form symmetric with respect to transposition $(1,2)$ of indices  (see Remark 1 in section \ref{secstructure}).  Let us choose the following invariants: 
$$
H=k^2=c^2, \ \ \ \omega_1=c(U_1+U_2), \ \ \ \omega_2=U_1-U_2.
$$
Around a point, where the Blaschke curvature does not vanish, the forms $\omega_1, \omega_2$ constitute a basis and one can differentiate with respect to this new basis. To keep the notation simple, we again denote this differentiation by sub-indices, avoiding confusion by explicitly introducing the derivatives.   
 Thus we define 
$g_1,g_2$, $g_{11},g_{12},g_{21},g_{22}$ as follows: 
$$
df=g_1\omega_1+g_2\omega_2,\ \ \  dg_1=g_{11}\omega_1+g_{12}\omega_2,\ \ \  dg_2=g_{21}\omega_1+g_{22}\omega_2.
$$
Similarly
$$
dH=H_1\omega_1+H_2\omega_2,\ \ \  dH_1=H_{11}\omega_1+H_{12}\omega_2.
$$
Now equations (\ref{structureEQk}) assume the form 
\begin{equation}\label{2pencilstructure}
d\omega_1=\frac{1}{2} \omega_1 \wedge \omega_2, \ \ \ d\omega_2=- \omega_1 \wedge \omega_2,\ \ \ H_2=-H,
\end{equation}
and equations (\ref{polymorphEQ}) read as
\begin{equation}\label{polymorphEQcontr}
g_{12}=g_2-\frac{1}{2}fg_1+\frac{1}{2}f^2+\frac{3}{4}f,\ \ \ g_{22}=\frac{H}{3}g_{11}+\left(\frac{H_1}{6}-H\right)g_1 +(f+1)g_2+\left(\frac{2H}{3}-\frac{H_1}{6}\right)f.
\end{equation}
\begin{lemma}\label{inZ2invariants}
Suppose that  functions $H,f$ and two 1-forms $\omega_1,\omega_2$ satisfy equations (\ref{2pencilstructure},\ref{polymorphEQcontr}) on some open set $V$, and that $H$ do not vanish on $V$.
Then the functions $f$, $k=\sqrt{H}$, $a=b=0$,  and 1-forms $U_1=\frac{1}{2k}(\omega_1+k\omega_2)$, $U_2=\frac{1}{2k}(\omega_1-k\omega_2)$ satisfy equations (\ref{structureEQk},\ref{polymorphEQ}).
\end{lemma}
{\it Proof:} One checks the claim by direct computation.
\hfill $\Box$\\

Now let us study the compatibility conditions of system (\ref{2pencilstructure},\ref{polymorphEQcontr}).
From $d(dH)=0$ and $d(df)=0$ we have
$$
H_{12}= H-\frac{H_1}{2}, \ \ \ g_{21}=g_{12}+g_2-\frac{g_1}{2}.
$$

Similarly, from $d(dg_1)=0$ and $d(dg_2)=0$ one obtains $g_{111}$ and $g_{112}$ via $f,g_1,g_2,g_{11},H,H_1$, where 
$$
dg_{11}=g_{111} \omega_1+g_{112} \omega_2.
$$
Now from $d(dg_{11})=0$ we get
$$
g_{11}=\frac{7(7f-4g_1)}{32}\frac{H_1}{H}-\frac{3f}{32}\frac{H_{11}}{H}+\frac{12f^3+36f^2+27f-296fH+288Hg_1-   24(4f+1)g_2}{64H}.
$$ 
With this expression for $g_{11}$, the condition $d(dg_{1})=0$ gives 
$$
\begin{array}{l}
H_{11}=\frac{1}{18f+48g_2}\{(640f^2-896fg_1+256g_1^2+3912f-1824g_1-576g_2)H+\\
\\
\ \ \ \ (96fg_1-96f^2-474f+24g_1+336g_2)H_1+\\
\\
\ \ \ \ 36f^3-480f^2g_2+108f^2-672fg_2-768g_2^2+81f-48g_2\}
\end{array}
$$
Substituting the above expression for $g_{11}$ into $d(dg_2)=0$,  we obtain an equation of the form:  
$$
T_{11}(g_1,g_2,f)H_1^2+T_{01}(g_1,g_2,f)H_1H+T_{00}(g_1,g_2,f)H^2+T_1(g_1,g_2,f)H_1+T_0(g_1,g_2,f)H=0,
$$
where the coefficients $T_{J}(g_1,g_2,f)$ are polynomial. 
Note that the equation is quadratic in $H_1,H$. For fixed $f,g_1,g_2$ we have a conic with one known point $(H,H_1)=(0,0)$. Parametrizing the conic by secants 
$$
H_1= ZH
$$
we express 
$$
H=\tilde{h}(g_1,g_2,f,Z),\ \ \ \ \ \  H_1=\tilde{h}_1(g_1,g_2,f,Z).
$$
With 
$$
dZ=Z_1\omega_1+Z_2\omega_2
$$ 
the form $Z\omega_1-\omega_2=\frac{dH}{H}$ must be closed, which gives 
$$
Z_2=\frac{1}{2}Z+1.
$$ 
Now substituting the obtained $H$ into
$$
dH=\tilde{h}_1(g_1,g_2,f,Z)\omega_1-\tilde{h}(g_1,g_2,f,Z)\omega_2
$$
one computes $Z_1=z_1(g_1,g_2,f,Z)$ as a rational function of  $g_1,g_2,f,Z$, equating the coefficients of $\omega_1$, and gets a polynomial equation $W(g_1,g_2,f,Z)=0$, equating the coefficients of $\omega_2$.   Observe that now  the differentials 
$dg_1,dg_2,dZ$ are expressed via $g_1,g_2,f,Z,\omega_1,\omega_2$: 
 \begin{equation}\label{Zg1g2}
\begin{array}{l}
dZ=z_1(g_1,g_2,f,Z)\omega_1+z_2(g_1,g_2,f,Z)\omega_2,\\
\\
dg_1=\hat{g}_{11}(g_1,g_2,f,Z)\omega_1+\hat{g}_{12}(g_1,g_2,f,Z)\omega_2,\\
\\
dg_2=\hat{g}_{21}(g_1,g_2,f,Z)\omega_1+\hat{g}_{22}(g_1,g_2,f,Z)\omega_2.
\end{array}
\end{equation}
Differentiate the last two equations of (\ref{Zg1g2}) and obtain two constraints $d(dg_i)=0$, 
involving $g_1,g_2,f,Z$; take the constraint numerators; factor the  resultant of these numerators with respect to $Z$; compute and also factor  such resultants of the numerators of $d(dg_i)=0$ with $W(g_1,g_2,f,Z)$. (In this computation we reduce the equation $d(dg_2)=0$ by a non-vanishing factor.) The obtained three sets of factors, considered without multiplicities,  intersect in a subset of five factors.  

 They are: $g_2$, $f-g_1$, $f-g_1+27/4$, $f^2+f+2g_2$, and  a factor $E(g_1,g_2,f)$ of degree 14. The four "simple" factors do not give non-hexagonal polymorphic webs, the compatibility analysis quickly implying  $f=0$. 
 
Thus we have to analyse the case $E(g_1,g_2,f)=0.$
Differentiating $E(g_1,g_2,f)$ and equating the coefficient of $\omega_1$ to zero, one obtains  $Z=z(g_1,g_2,f)$ as a rational (and rather involved) function. Comparing the coefficient of $\omega_2$ in $dz(g_1,g_2,f)$  with $Z_2=\frac{1}{2}Z+1$, we get one more polynomial equation $\tilde{E}(g_1,g_2,f)=0$ of degree 77.
 
\begin{theorem}
There is no non-hexagonal polymorphic 3-web with 2 pencils of lines.
\end{theorem}
{\it Proof:} The polynomials $E,\tilde{E}$ are irreducible over $\mathbb Q$, therefore they define some curve in 3-dimensional space. Then all the projective invariants of the web are parametrized by points on this curve, and the signature set is not "larger" than one-dimensional. By Theorem \ref{1dim}, the web admits at least one infinitesimal projective symmetry. But there is no  polymorphic 3-web admitting infinitesimal symmetry (see \cite{Agr}). \hfill $\Box$

\section*{Acknowledgement}

This research was supported by grants \#2014/17812-0 and  \#2017/02954-2 of S\~ao Paulo Research Foundation (FAPESP).

\section*{Appendix}
Here we present the details of the computations, mentioned in the proofs.
\subsection{Compatibility conditions in the proof of Lemma \ref{aligned}}
From equations (\ref{Euler}) one gets all mixed derivatives in terms of derivatives only in $y$. 
Resolving the constraint (\ref{constraintLine}) for $R_y$, one obtains
$$
R_y=[(R-Q)P_y+(P-R)Q_y]/(P-Q).
$$
 Now comparing $R_{xy}$, obtained from (\ref{Euler}), with 
$$
R_{xy}=D_x([(R-Q)P_y+(P-R)Q_y]/(P-Q))
$$ yields $Q_{yy}=P_{yy}$. (Here $D_x=\partial_x+P_x\partial_P+Q_x\partial_Q+R_x\partial_R + ...$ is the operator of total derivative with respect to $x$.)
Thus $Q_{xyy}=P_{xyy}$, which implies $P_{yyy}= 3P_{yy}(Py-Qy)/(Q-P)$. Finally,
$$
P_{xyyy}= 3D_x(P_{yy}(Py-Qy)/(Q-P))
$$
gives $P_{yy}=0$ and therefore $Q_{yy}=0$. Hence  $R_{yy}=0$ due to the permutation symmetry.

\subsection{Computing compatibility conditions for the proof of Theorem \ref{analyticdim}}
The length of the expressions, involved in the computation of compatibility conditions, grows very quickly. Therefore we explain here, step by step, the computation scheme without giving explicit formulas. This scheme was implemented on a notebook with 16GB of memory. The software used was Maple 18.

\noindent First, we write $k_1$ and $k_2$ via $h$, using the last equation of (\ref{structureEQk}), and find $f_{11},f_{22}$ from equations (\ref{polymorphEQ}). Then the equations $d(df_1)=d(df_2)=0$ give $L_1$ and $L_2$. 

\noindent The condition $d(dk)=0$ gives $h_1,h_2$ via $r:=\frac{h_2-h_1}{2}$.

\noindent Finding $r_1$ from $d(dh)=0$ and substituting it into $d(dL)=0$, we obtain $r$ and, consequently, $h_1$ and $h_2$. 

\noindent The condition $d(dh)=0$ (Note that, with $h_1$ and $h_2$ found, we have to differentiate $dh$ again!) permits to express $a_{222},b_{111}$ via $m:=\frac{b_{111}-a_{222}}{2}$.

\noindent Equations $d(da_{22})=d(db_{11})=0$ give $m_1$ and $m_2$. 

\noindent Now $d(dm)=0$ gives $m$.  Finally, differentiating again the expressions for $da_{22}$ and $db_{11}$ we get two independent polynomial equations for 12 invariants $a,a_2,a_{22},b,b_1,b_{11},k,h,f,f_1,f_2,L$ from   $d(da_{22})=d(db_{11})=0$.


\begin{thebibliography}{99}


\addcontentsline{toc}{section}{References}

\leftmargin=14cm
\labelwidth=1.3cm


\bibitem{Agr} Agafonov S.I., Gronwall's conjecture for 3-webs with  infinitesimal symmetries, (2014) arXiv: 1411.0874 [math.DG],  to appear in Comm. Anal. Geom. 28 (2020), no 5.

\bibitem{AGL} Akivis, M.A., Goldberg, V.V., Lychagin, V.V., Linearizability of $d$-webs, $d > 4$, on two-dimensional manifolds,
{\it Selecta Math.} 10(4) (2004), 431--451.

\bibitem{BB} Blaschke, W., Bol, G.,
{\it Geometrie der Gewebe, Topologische Fragen der
Differentialgeometrie} J. Springer, Berlin, 1938.

\bibitem{Be} Blaschke, W., {\it Einf\"uhrung in die Geometrie der
Waben}, Birkh\"auser Verlag, Basel und Stuttgart, 1955.

\bibitem{Bk} Bol, G., Geradlinige Kurvengewebe.
"Topologische fragen der differentialgeometrie 31."
{\it Abh. Math. Sem. Univ. Hamburg}, 8, (1931), no. 1, 264--270.

\bibitem{Bg} Bol, G.,
Ueber Geradengewebe.
"Topologische Fragen der Differentialgeometrie (65)."
{\it Ann. Mat. Pura Appl.}, 17 (1938), no. 1, 45--58.

\bibitem{Bp} Bor\r{u}vka, O., Sur les correspondances analytiques entre deux plans projectifs II, {\it Univ. Mazaryk \u{C}.},
85, (1938)
22--24.

\bibitem{Ce} Cartan, \`E., {\it Les syst\`emes diff\'erentiels ext\'erieurs et leurs applications g\'eom\'etriques.} (French) 
 Hermann et Cie., Paris, 1945. 

\bibitem{GLb} Goldberg, V.V., Lychagin, V.V., On the Blaschke conjecture for 3-webs. {\it J. Geom. Anal.} 16 (2006), no. 1, 69--115.

\bibitem{GSg} Graf, H., Sauer. R.,  \"Uber dreifache
Geradensysteme in der Ebene, welche Dreiecksnetze bilden, {\it
Sitzungsb. Math.-Naturw. Abt.} (1924), 119--156.

\bibitem{Gn} Gronwall, T.H., Sur les \'equations entre trois variables repr\'esentables par
les nomogrammes \`a points align\'e, {\it J. de Liouville}, 8, (1912), 59--102.

\bibitem{GMS} Grifone, J., Muzsnay, Z., Saab, J., On the linearizability of 3-webs, Proceedings of the Third World Congress
of Nonlinear Analysis, Part 4 (Catania, 2000), {\it Nonlinear Anal,} 47(4) (2001), 2643--2654.

\bibitem{Hl} H\'enaut, A.,
Sur la lin\'earisation des tissus de $\mathbb C^2$.
{\it Topology 32} (1993), no. 3, 531--542.

\bibitem{Lg}  Liouville R., Sur une classe d'\'equations diff\'erentielles, parmi lesquelles, en particulier,
toutes celles des lignes g\'eod\'esiques se trouvent comprises, {\it Comptes rendus
hebdomadaires des s\'eances de l'Acad\'emie des sciences} 105 (1887) 1062--1064.

\bibitem{Og} Ovsiannikov, L.V.,  {\it Group analysis of differential equations.}, New York-London, 1982.

\bibitem{Pl} Pirio, L., Sur la lin\'earisation des tissus. {\it Enseign. Math.} (2) 55 (2009), no. 3--4, 285--328.

\bibitem{Sp} Sasaki, T., {\it Projective Differential Geometry and Linear Homogeneous Differential Equations}, Rokko Lectures in Math., 5. Kobe University, 1999.

\bibitem{Sc} Sharpe, R.W., {\it Differential geometry. Cartan's generalization of Klein's Erlangen program.} Graduate Texts in Mathematics, 166. Springer-Verlag, New York, 1997.

\bibitem{Su} Smirnov, S.V., On certain problems of uniqueness in the theory of webs. (Russian) {\it Vol. Mat. Sb.}, 2 (1964),
128--135.

\bibitem{Sa} Smirnov, S.V., Uniqueness of a nomogram of aligned points with one rectilinear scale. (Russian)
{\it Sibirsk. Mat. $\check{ Z}.$,} 5 (1964), 910--922.

\bibitem{Vn} Vaona, G., Sur teorema fondamentale della nomografia, {\it Boll. Un. Mat. ltal.,} (3) 16 (1961), 258--263.

\bibitem{Vi} Vessiot, E., Sur l'int\'egration des syst\`emes diff\'erentiels qui admettent des groupes continus de transformations. {\it Acta Math.} 28 (1904), no. 1, 307--349.

\bibitem{Wg} Wang, J.S., On the Gronwall conjecture. {\it J. Geom. Anal.} 22 (2012), no. 1, 38--73.

\bibitem{Wpc} Wilczynski, E.J., {\it Projective differential geometry of curves and ruled surfaces,} (German) Teubner (1906).

\end{thebibliography}
\end{document}